\documentclass[
12pt,reqno]{amsart}
\usepackage{verbatim}
\usepackage{amssymb}

\usepackage{enumerate}
\usepackage[active]{srcltx}
\numberwithin{equation}{section}

\usepackage{t1enc}
\usepackage[utf8x]{inputenc}

\newtheorem{theorem}{Theorem}[section]
\newtheorem{lemma}[theorem]{Lemma}
\newtheorem{corollary}[theorem]{Corollary}

\newtheorem{proposition}[theorem]{Proposition}

\theoremstyle{definition}
\newtheorem{definition}[theorem]{Definition}

\theoremstyle{remark}

\newcommand{\itemprefix}{}
\makeatletter
\newcommand{\myitem}{%
\item\protected@edef\@currentlabel{\itemprefix\theenumi}%
}
\makeatother

\newcommand{\setm}{\setminus}
\newcommand{\empt}{\emptyset}
\newcommand{\subs}{\subset}
\newcommand{\dom}{\operatorname{dom}}

\newcommand{\<}{\left\langle}
\renewcommand{\>}{\right\rangle}

\author[I. Juh\'asz]{Istv\'an Juh\'asz}
\address      { Alfr\'ed Rényi Institute of Mathematics%
}
\email{juhasz@renyi.hu}

\author[S. Shelah]{Saharon Shelah}
\address
      { Institute of Mathematics\\Hebrew University, Jerusalem  }

\author[L. Soukup]{Lajos Soukup}
\thanks
  {
   }
\address
      { Alfr{\'e}d R{\'e}nyi Institute of Mathematics
}
\email{soukup@renyi.hu}

\author[Z. Szentmikl\'ossy]{Zolt\'an Szentmikl\'ossy}
\address{E\"otv\"os University of Budapest}
\email{szentmiklossyz@gmail.com}
\thanks{The first, third and fourth authors were supported
by NKFIH grant no. K129211. Research partially supported by the Israel Science Foundation (ISF)
grant no: 1838/19.
 Research partially supported by NSF grant no: DMS 1833363.
 Publication number Sh1213.
}


\title{Large strongly anti-Urysohn spaces exist}

\begin{document}

\begin{abstract}
As defined in \cite{JSSz}, a Hausdorff space is {\em strongly anti-Urysohn} (in short: SAU) if it has at least two non-isolated
points and any two {\em infinite} closed subsets of it intersect. Our main result answers the two main questions of \cite{JSSz} by
providing a ZFC construction of a locally countable SAU  space of cardinality $2^{\mathfrak{c}}$.
The construction hinges on the existence of $2^{\mathfrak{c}}$ weak P-points in $\omega^*$,
a very deep result of Ken Kunen.

It remains open if SAU  spaces of cardinality $> 2^{\mathfrak{c}}$ could exist, while it was shown in \cite{JSSz} that
$2^{2^{\mathfrak{c}}}$ is an upper bound. Also, we do not know if {\em crowded} SAU  spaces,
i.e. ones without any isolated points, exist in ZFC but
we obtained the following consistency results concerning such spaces.

\smallskip

\begin{enumerate}
\item
It is consistent that $\mathfrak{c}$  is  as large
as you wish  and
there is a locally countable and crowded SAU space
of cardinality $\mathfrak{c}^+$.

\smallskip

\item It is consistent that both $\mathfrak{c}$ and $2^\mathfrak{c}$ are as large
as you wish  and
there is a crowded SAU space
of cardinality $2^\mathfrak{c}$.

\smallskip

\item
 For any uncountable cardinal ${\kappa}$ the following statements are equivalent:
\begin{enumerate}[(i)]
      \item ${\kappa}=cof({[{\kappa}]}^{{\omega}},\subseteq)$.
      \item
      There is a locally countable and
      crowded SAU space of
      size ${\kappa}$ in the generic extension obtained
      by adding $\kappa$ Cohen reals.
      \item There is a locally countable and
      countably compact $T_1$-space of
      size ${\kappa}$ in some CCC generic extension.
      \end{enumerate}
\end{enumerate}

\end{abstract}

\dedicatory{Dedicated to the memory of our old friend Ken Kunen}

\maketitle

\section{Introduction}

Anti-Urysohn (AU) and strongly anti-Urysohn (SAU) spaces were introduced and studied in
\cite{JSSz}. An AU space is a Hausdorff space in which any two non-empty {\em regular} closed sets intersect
and a SAU space is a Hausdorff space that has at least two non-isolated
points and in which any two {\em infinite} closed sets intersect.
Note that a non-singleton AU space has no isolated points, i.e. is crowded and a crowded SAU space is AU,
this explains the terminology. Also, the requirement of having at least two non-isolated
points is needed to exclude the trivial case of one-point compactifications of discrete spaces.

All relevant questions concerning AU spaces were settled in \cite{JSSz}, in particular it was
shown that for every infinite cardinal $\kappa$ there is an AU space of cardinality $\kappa$, but
only inconclusive partial results were proved for SAU spaces. For instance, we could only construct consistent
examples of SAU spaces, moreover all of them were of size $\le \mathfrak{c}$, while only $2^{2^\mathfrak{c}}$
was established as an upper bound for their cardinality.

This, of course, naturally led to the following two questions raised in \cite{JSSz}:

\begin{enumerate}
\item Is there a SAU space in ZFC?

\smallskip

\item Can there be a SAU space of cardinality greater than $\mathfrak{c}$?
\end{enumerate}
Our main result answers affirmatively both of these questions, namely we shall present a ZFC example
of a locally countable SAU space of cardinality $2^\mathfrak{c}$. It is easy to see that $2^\mathfrak{c}$
is an upper bound for the sizes of locally countable SAU spaces, so this result is sharp. However, it
remains an open question if $2^\mathfrak{c}$
is an upper bound for the sizes of all SAU spaces. It was proved in \cite{JSSz} that
$2^{2^{\mathfrak{c}}}$ is such an upper bound.

It was proved in   \cite[Theorem 3.7]{JSSz} that $\mathfrak{r} = \mathfrak{c}$ implies the existence of a locally countable
and {\em crowded} SAU space  of size $\le \mathfrak{c}$, moreover several other forcing constructions yielded
{\em crowded} SAU spaces. This then led to the question if the existence of a SAU space is equivalent to
the existence of crowded ones.

By Theorem \ref{tm:SAU} this question can now be reformulated as follows: Is there a crowded SAU space in ZFC?
Now, the method of construction of Theorem \ref{tm:SAU} yields a SAU space that is right-separated, i.e. scattered,
hence it may not help to answer this question that we could not answer. However, we could get several
partial results about it. We could prove the consistency of the existence of a locally countable
and {\em crowded} SAU space  of size $\le \mathfrak{c}^+$, moreover we proved that the equality
${\kappa}=cof({[{\kappa}]}^{{\omega}},\subseteq)$ is equivalent to the existence of a locally countable
and {\em crowded} SAU space in the generic extension obtained by adding $\kappa$ Cohen reals.

\section{A large SAU space in ZFC}

The aim of this section is to present a construction that yields, in ZFC,
a locally countable SAU space of cardinality $2^\mathfrak{c}$.We actually
present a general construction of right-separated spaces that uses functions
with values that are free filters. We start by fixing some notation and
terminology.

For any infinite set $S$ we write $\Phi(S)$ to denote the collection of all
free filters on $S$. As is customary, $S^* \subs \Phi(S)$ denotes the family
of all free ultrafilters on $S$.

\begin{definition}\label{df:filt}
Let $\kappa$ be an uncountable cardinal. We call a function $\varphi$ with domain $\kappa \setm \omega$ a
{\em nice filter assignment} for $\kappa$ if for all $\alpha \in \kappa \setm \omega$ we have
$\varphi(\alpha) \in \Phi(S_\alpha)$ for some infinite subset $S_\alpha$ of $\alpha$.
We shall denote by $\mathfrak{F}(\kappa)$ the family of all nice filter assignments for $\kappa$.
\end{definition}
Note that $S_\alpha$ is simply the maximal member of $\varphi(\alpha)$.

\smallskip

Next we shall show how a nice filter assignment for $\kappa$ defines a natural topology on $\kappa$.

\begin{definition}\label{df:ftop}
Let $\varphi$ be a nice filter assignment for $\kappa > \omega$. Then the topology $\tau_\varphi$ on $\kappa$
is defined by the formula $$\tau_\varphi = \{G \subs \kappa : \forall\,\alpha \in G \setm \omega\, \big(G \cap S_\alpha \in \varphi(\alpha)\big)  \}.$$
\end{definition}

It is left to the reader to check that $\tau_\varphi$ is a $T_1$ topology on $\kappa$, this is where
considering only free filters is essential.

Since SAU spaces are Hausdorff by definition, we shall need a condition that will imply in case of
a nice filter assignment $\varphi$ for $\kappa$ that $\tau_\varphi$ is Hausdoff.
To formulate this condition, we shall use the following terminology that was introduced in \cite{HJ}.

\begin{definition}\label{df:drep}
The indexed family of filters $\{F_i : i \in I\}$ is called {\em disjointly representable} if there
are sets $\{A_i \in F_i : i \in I\}$ such that $A_i \cap A_j = \emptyset$ whenever
$i,j \in I$ and $i \ne j$.
\end{definition}

\begin{lemma}\label{lm:T2}
Let $\varphi$ be a nice filter assignment for $\kappa > \omega$ such that
\begin{enumerate}[(i)]
\item $S_\alpha$ is countable for all $\alpha \in \kappa \setm \omega$,

\smallskip

\item $\{\varphi(i) : i \in I\}$ is disjointly representable for all countable $I \subs \kappa \setm \omega$.
\end{enumerate}
Then $\tau_\varphi$ is a Hausdoff topology.
\end{lemma}

\begin{proof}
Note first that all $n \in \omega$ are isolated in $\tau_\varphi$. So, it suffices to show that
distinct $\alpha, \beta \in \kappa \setm \omega$ have disjoint neighborhoods.

To see that, we set $U_0 = \{\alpha\}$ and $V_0 = \{\beta\}$ and then define by recursion countable sets $U_n$ and $V_n$
with $U_n \cap V_n = \emptyset$ as follows.

Given  $U_n$ and $V_n$, we apply (ii) to the index set $I_n = U_n \cup V_n \setm \omega$ to obtain
pairwise disjoint sets $A_i \in \varphi(i)$ for all $i \in I_n$. Then we set $U_{n+1} = U_n \cup \bigcup \{A_i : i \in U_n\}$
and similarly $V_{n+1} = V_n \cup \bigcup \{A_i : i \in V_n\}$. Clearly, then both $U_{n+1}$ and $V_{n+1}$ are countable by (i),
moreover we also have $U_{n+1} \cap V_{n+1} = \emptyset$. Now, it is obvious that $U = \bigcup \{U_n : n< \omega\} \in \tau_\varphi$ and
$V = \bigcup \{V_n : n< \omega\} \in \tau_\varphi$ are disjoint open neighborhoods of
$\alpha$ and $\beta$, completing the proof.
\end{proof}

In view of this result it is natural to look for conditions that imply disjoint 
representability of certain
countable families of filters. Here is a very simple such condition about families
 of ultrafilters on $\omega$.

\begin{proposition}\label{pr:dis}
A countable subfamily of $\omega^*$ is disjointly representable iff it is a discrete subspace of $\omega^*$, considered as the
remainder of $\beta\omega$.
\end{proposition}

\begin{proof}
Indeed, this follows from the fact that in a regular space the points in a countable discrete subspace have
pairwise disjoint neighborhoods, moreover for members of $\omega^*$ this means that they have
pairwise almost disjoint elements.
\end{proof}

Contrary to this proposition, the following result needed in our construction of SAU spaces, is highly
non-trivial. But, luckily for us, it is an immediate consequence of a deep result of Kunen in \cite{K}.
We recall that an ultrafilter $u \in \omega^*$ is a {\em weak P-point} if $u$ is not in the closure
of any {\em countable} subset of its complement.

\begin{theorem}\label{tm:Ku}
The family $\mathcal{U}$ of all weak P-point ultrafilters in $\omega^*$ has cardinality $2^\mathfrak{c}$,
moreover
all countable subfamilies of $\mathcal{U}$ are disjointly representable.
\end{theorem}

\begin{proof}
We have $|\mathcal{U}| = 2^\mathfrak{c}$
by \cite{K} and it is obvious that all countable subsets of $\mathcal{U}$ are
discrete in $\omega^*$, hence disjointly representable by Proposition \ref{pr:dis}.
\end{proof}

We are now ready to present our main result.

\begin{theorem}\label{tm:SAU}
If $\kappa = \kappa^\omega \le 2^\mathfrak{c}$ then there is a locally countable SAU space of cardinality $\kappa$.
\end{theorem}

\begin{proof}
To start with, we fix using Theorem \ref{tm:Ku} for every countably infinite set $S \subs \kappa$
the family $\mathcal{U}(S)$ of size $2^\mathfrak{c}$ of all weak P-point ultrafilters in $S^*$.  Then
all countable subfamilies of $\mathcal{U}(S)$ are disjointly representable. Also, if $T \in [S]^\omega$
then $$\mathcal{U}(T) = \{u\upharpoonright T : u \in \mathcal{U}(S) \text{ and } T \in u\},$$
where $u\upharpoonright T = \{A \cap T : A \in u\}$. In other words, $u \in T^*$ belongs to $\mathcal{U}(T)$
iff the ultrafilter $\widehat{u} \in S^*$ generated by $u$ belongs to $\mathcal{U}(S)$.
Indeed, this is because $T^*$ is clopen in $S^*$.

Next, using $\kappa = \kappa^\omega$ we fix an enumeration $\{\langle A_\alpha,B_\alpha \rangle : \omega \le \alpha < \kappa\}$
of $[\kappa]^\omega \times [\kappa]^\omega$
such that $S_\alpha = A_\alpha \cup B_\alpha \subs \alpha$ for all $\alpha \in \kappa \setm \omega$.
We then pick, by transfinite recursion on $\alpha \in \kappa \setm \omega$, weak P-points $u_\alpha \in \mathcal{U}(S_\alpha)$
and $v_\alpha \in \mathcal{U}(S_\alpha)$ such that $A_\alpha \in u_\alpha$ and $B_\alpha \in v_\alpha$ as follows.

Assume that $\alpha \in \kappa \setm \omega$ and $u_\beta,v_\beta$ have been chosen for all $\omega \le \beta < \alpha$
and then let us put $$I = \{\beta \in \alpha \setm \omega : A_\beta \cap S_\alpha \in u_\beta\} \text{ and }
J = \{\beta \in \alpha \setm \omega : B_\beta \cap S_\alpha \in v_\beta\}.$$
For each $\beta \in I$ then $u_\beta \upharpoonright A_\alpha$ generates a weak P-point ultrafilter $\widehat{u}_\beta \in \mathcal{U}(S_\alpha)$
and, similarly, for each $\beta \in J$, $\,v_\beta \upharpoonright B_\alpha$ generates a weak P-point ultrafilter $\widehat{v}_\beta \in \mathcal{U}(S_\alpha)$.

But we have on one hand $|I \cup J| < \kappa \le 2^\mathfrak{c}$ and, on the other hand, $ |\mathcal{U}(S_\alpha)| = |\mathcal{U}(A_\alpha)| = |\mathcal{U}(B_\alpha)| = 2^\mathfrak{c}$,
so we may clearly choose distinct $u_\alpha,v_\alpha \in \mathcal{U}(S_\alpha) \setm \{\widehat{u_\beta} : \beta \in I\} \cup \{\widehat{v_\beta} : \beta \in J\}$
such that $A_\alpha \in u_\alpha$ and $B_\alpha \in v_\alpha$.

After having completed the recursion, we let $$\varphi(\alpha) = \{U \cup V : U \in u_\alpha \text{ and } V \in v_\alpha\}$$
for each $\alpha \in \kappa \setm \omega$. Clearly, then $\varphi(\alpha) \in \Phi(S_\alpha), hence$ $\varphi$ is a nice filter assignment for $\kappa$.
It is also clear from the definitions that each $\alpha \in \kappa \setm \omega$ is a common $\tau_\varphi$-accumulation point
of both $A_\alpha$ and $B_\alpha$, hence $\tau_\varphi$ turns out to be a SAU topology on $\kappa$ if we can prove that it is Hausdorff.

To see that, it suffices to show that for any $I \in [\kappa \setm \omega]^\omega$ the family $\{u_\alpha,v_\alpha : \alpha \in I\}$
is disjointly representable. Indeed, if $$\{U_\alpha : \alpha \in  I\} \cup \{V_\alpha : \alpha \in  I\}$$ are pairwise disjoint sets with
$U_\alpha \in u_\alpha$ and $V_\alpha \in v_\alpha$  then $\{U_\alpha \cup V_\alpha : \alpha \in I\}$ are pairwise disjoint as well.
But this means that $\{\varphi(\alpha) : \alpha \in I\}$ is disjointly representable for all countable $I \subs \kappa \setm \omega$,
hence $\tau_\varphi$ is Hausdorff by Lemma \ref{lm:T2}.

So, consider $I \in [\kappa \setm \omega]^\omega$ and put $S = \bigcup \{S_\alpha\ : \alpha \in I\}$. For each $\alpha \in I$
then $u_\alpha$ generates a weak P-point $\widehat{u}_\alpha \in \mathcal{U}(S)$ and similarly $v_\alpha$ generates
$\widehat{v}_\alpha \in \mathcal{U}(S)$, moreover by our recursive construction they are all distinct.
Consequently, by Theorem \ref{tm:Ku} the family $\{\widehat{u}_\alpha, \widehat{v}_\alpha : \alpha \in I\}$ is
disjointly representable, hence so is $\{u_\alpha,v_\alpha : \alpha \in I\}$, completing the proof.
\end{proof}

It is worth to mention that for $\kappa \ge \mathfrak{c}$ the condition $\kappa = \kappa^\omega \le 2^\mathfrak{c}$ in Theorem \ref{tm:SAU} is
actually necessary to have a locally countable SAU space of cardinality $\kappa$. Since SAU spaces are countably compact, this follows immediately from
the following result.

\begin{theorem}\label{tm:kw}
If $X$ is a locally countable and countably compact $T_1$-space of cardinality $\kappa > \mathfrak{c}$ then $\kappa = \kappa^\omega$.
\end{theorem}

\begin{proof}
Assume, on the contrary that $\kappa < \kappa^\omega$ and let $\lambda$ be the smallest cardinal such that $\lambda^\omega > \kappa$.
It is well-known that then $cf(\lambda) = \omega$. But then by a classical
 result of Tarski in \cite{T}, there is an almost disjoint
family $\mathcal{A} \subs [\lambda]^\omega$ with $|\mathcal{A}| = \lambda^\omega > \kappa$.

Since $\lambda \le \kappa$, we may assume without any loss of generality that $\lambda \subs X$, hence every $A \in \mathcal{A}$
has an accumulation point $x_A \in X$. But $|\mathcal{A}| > \kappa$ then implies the existence of some $\mathcal{B} \subs \mathcal{A}$
and $x \in X$ such that $|\mathcal{B}| > \kappa$ and $x_A = x$ for all $A \in \mathcal{B}$. Let $U$ be a countable neighbourhood of $x$,
then $A \cap U$ is infinite for all $A \in \mathcal{B}$ which is impossible because $\mathcal{B}$ is almost disjoint
with $|\mathcal{B}| > \kappa > \mathfrak{c}$.
\end{proof}

\begin{corollary}
For $\mathfrak{c} \le \kappa \le 2^\mathfrak{c}$ there is a locally countable SAU space of cardinality $\kappa$
iff $\kappa = \kappa^\omega$.
\end{corollary}

As was mentioned in the introduction, any locally countable SAU space $X$ has cardinality $\le 2^\mathfrak{c}$.
Indeed, this is because it does have an infinite closed subset $F$ of cardinality $\le 2^\mathfrak{c}$, namely
the closure of any subset of size $\omega$. But local countability then implies that $F$ is covered by an open
set $U$ with $|U| = |F| \le 2^\mathfrak{c}$ and the SAU property implies that $X \setm U$ is finite.

Finally, we mention the following accidental consequence of Theorem \ref{tm:SAU}. In this, $F(X)$ is the free set number and
$X_\delta$ denotes the $G_\delta$-modification of the space $X$, see \cite{JSSz1}.

\begin{corollary}\label{co:F}
The locally countable SAU space $X$ of cardinality $2^\mathfrak{c}$ of Theorem \ref{tm:SAU} is an
example of a Hausdorff space with $F(X) = \omega$
and $F(X_\delta) = 2^\mathfrak{c}$.
\end{corollary}

Indeed, $F(X) = \omega$ because every free sequence in a SAU space has order type $< \omega + \omega$, while $X_\delta$ is discrete.
It was shown in \cite{JSSz1} that for any Hausdorff space $X$ with $F(X) = \omega$ we have $F(X_\delta) \le 2^{2^\mathfrak{c}}$
and we do not know if this upper bound could be replaced by $2^\mathfrak{c}$, i.e. if Corollary \ref{co:F} is sharp.
It is curious that the same upper bound $2^{2^\mathfrak{c}}$ is known for the size of any SAU space and the same problem
arises if this could be improved to $2^\mathfrak{c}$.

\bigskip

\section{Forcing "large" crowded SAU spaces}

All the {\em consistent} examples of SAU spaces constructed in \cite{JSSz} were crowded but
of cardinality $\le \mathfrak{c}$. As mentioned above, we do not know if {\em crowded} SAU spaces
exist in ZFC but the aim of this section is to produce consistent examples of crowded SAU spaces
of size $> \mathfrak{c}$.

By Theorem 3.7 of \cite{JSSz}, the assumption $\mathfrak{r} = \mathfrak{c}$ implies the existence
of a locally countable and crowded SAU space of cardinality $\mathfrak{c}$. Our next result says that, under
the same assumption $\mathfrak{r} = \mathfrak{c}$ together with $\mathfrak{c} = 2^{<\mathfrak{c}}$,
a forcing construction yields a generic extension of the ground model
in which there is a locally countable and crowded SAU space of cardinality $\mathfrak{c}^+$.

\begin{theorem}\label{tm:r=c}
If $\mathfrak{r} = \mathfrak{c} = 2^{<\mathfrak{c}}$ then we have a $\mathfrak{c}$-closed and
$\mathfrak{c}^+-CC$ notion of forcing $\mathbb{P}$ such that, in the generic extension $V^\mathbb{P}$, there is
a locally countable and crowded SAU space of cardinality $\mathfrak{c}^+$.
\end{theorem}

\begin{proof}
Our aim is to obtain a function $U : \mathfrak{c}^+ \to [\mathfrak{c}^+]^\omega$ in $V^\mathbb{P}$ such that $\alpha \in U(\alpha)$
for each $\alpha \in \mathfrak{c}^+$ and $\{U(\alpha) : \alpha \in \mathfrak{c}^+\}$ generates a SAU topology on $\mathfrak{c}^+$.
Our conditions then will be approximations to $U$ of size $< \mathfrak{c}$ with some ``side conditions" that will ensure that
any two infinite subsets of $\mathfrak{c}^+$ have a common accumulation point. 
Hausdorffness of $\tau$ will follow
from the assumption $\mathfrak{r} = \mathfrak{c}$ and genericity.

Now we define the notion of forcing $\mathbb{P} = \<P, \le \>$. The elements of $P$ will be pairs of the form
$p = \<U_p,\mathcal{C}_p\>$, where $U_p$ is a function with
domain $A_p \in [\mathfrak{c}^+]^{< \mathfrak{c}}$ and values taken in $[A_p]^\omega$ such that $\omega \subset A_p$ and
$U_p(n) = \omega$ for all $n \in \omega$, moreover $\alpha \in U(\alpha) \subs \alpha+1$ if $\alpha \in A_p \setminus \omega$.

Then $\{U(\alpha) : \alpha \in A_p\}$ as a subbase generates a topology $\tau_p$ on $A_p$ that is required to be crowded,
i.e. we assume that all non-empty members of $\tau_p$ are infinite. To handle this, we define
$$B_{p,I} = \bigcap \{U_p(\alpha) : \alpha \in I\}$$ for any finite subset $I$ of $A_p$.
Then $$\mathcal{B}_p = \{B_{p,I} : I \in [A_p]^{< \omega}\} \setminus \{\emptyset\}$$ is a base for $\tau_p$, hence our assumption
just means that all members of $\mathcal{B}_p$ are infinite.

For any $x \in A_p$ we shall denote by $ac_p(x)$ the family of all sets $C \in [A_p]^\omega$ such that $x$ is a
$\tau_p$-accumulation point of $C$, i.e. every $\tau_p$-neighborhood of $x$ has infinite intersection with $C$.
Now, the second co-ordinate $\mathcal{C}_p$ of the condition $p$ is also a function with domain $A_p$ but
such that $\mathcal{C}_p(x) \in [ac_p(x)]^{< \mathfrak{c}}$ for $x \in A_p$.

Next we define the partial order $\le$ on $P$ by the following stipulations:

\begin{definition}\label{df:leq1}
For $p,q \in P$ we have $p \leq q$, i.e. $p$ is a stronger condition than $q$ iff
\begin{enumerate}[(a)]
\item $U_p \supset U_q$ and

\smallskip

\item $\mathcal{C}_p(x) \supset \mathcal{C}_q(x)$ for all $x \in A_q$.
\end{enumerate}
\end{definition}

We may reformulate item (a) above as follows: $A_p \supset A_q$ and $U_p(x) = U_q(x)$ for all $x \in A_q$.
Note that this implies $B_{p,I} = B_{q,I}$ for all $I \in [A_q]^{< \omega}$.
Also, it should be noted that in item (b) it is implicit that every $C \in  \mathcal{C}_q(x)$ accumulates
to $x$ in the topology $\tau_p$ as well.

We next present several lemmas concerning the forcing $\mathbb{P}$ which together will yield the desired
function $U : \mathfrak{c}^+ \to [\mathfrak{c}^+]^\omega$  in $V^\mathbb{P}$.

\begin{lemma}\label{lm:c+cl1}
The forcing $\mathbb{P}$ is $\mathfrak{c}$-closed.
\end{lemma}

\begin{proof}[Proof of \ref{lm:c+cl1}]
Assume that $\<p_\xi : \xi < \varrho\>$ is a decreasing $\varrho$-sequence in $\mathbb{P}$, where
$p_\xi = \< U_\xi, \mathcal{C}_\xi \>$ for $\xi < \varrho$ and $\varrho < \mathfrak{c}$  is an infinite regular cardinal.
(To enhance readability, we use the $\xi$'s instead of the $p_\xi$'s
as indeces.) We may then define a lower bound $p = \<U_p,\mathcal{C}_p\>$ for the $p_\xi$'s as follows.

\begin{enumerate}
\item  $U_p = \bigcup_{\xi < \varrho} U_\xi$, hence $A_p = \bigcup_{\xi < \varrho} A_\xi$, and

\smallskip

\item $\mathcal{C}_p(x) = \bigcup \{\mathcal{C}_\xi(x) : x \in A_\xi \}$ for any $x \in A_p$.
\end{enumerate}

To see that $p \in P$, first note that $\mathfrak{c} = 2^{<\mathfrak{c}}$ implies that $\mathfrak{c}$
is regular, hence we have both $|A_p| < \mathfrak{c}$ and $|\mathcal{C}_p(x)| < \mathfrak{c}$ for any $x \in A_p$.
That $\tau_p$ is crowded follows from the fact that $\mathcal{B}_p = \bigcup_{\xi < \varrho} \mathcal{B}_\xi$.
It remains to check that for any $C \in \mathcal{C}_\xi(x)$ we have $C \in ac_p(x)$, and that is
clear because then $C \in ac_\eta(x)$ for all $\eta \in \varrho \setminus \xi$ and every $B \in \mathcal{B}_p$
eventually belongs to $\mathcal{B}_\eta$ as well.
\end{proof}

Our next lemma is an amalgamation result that will imply the $\mathfrak{c}^+-CC$-property of $\mathbb{P}$.

\begin{lemma}\label{lm:amalg1}
Assume that $p,q \in P$ are such that $A_p \cap A_q < A_p \Delta A_q$, $\,otp(A_p) = otp(A_q)$, and for
the unique order isomorphism $\pi : A_p \to A_q$ between them we have, for all $\alpha \in A_p$, that

\begin{enumerate}[(i)]
\item $U_q(\pi(\alpha)) = U_p(\alpha)$, and

\smallskip

\item $\mathcal{C}_q(\pi(\alpha)) = \{\pi[C] : C \in \mathcal{C}_p(\alpha)\}$.
\end{enumerate}
Then $p$ and $q$ are compatible in $\mathbb{P}$.
\end{lemma}

\begin{proof}[Proof of \ref{lm:amalg1}]
Let us note first that, as $\pi$ is the identity on $A_p \cap A_q < A_p \Delta A_q$, we have
$U_p(\alpha) = U_q(\alpha)$ for all $\alpha \in A_p \cap A_q$ by (i), hence $U_r = U_p \cup U_q$
is a well-defined map on $A_r = A_p \cup A_q$. We claim that if we set $\mathcal{C}_r(x) = \mathcal{C}_p(x) \cup \mathcal{C}_q(x)$
for $x \in A_p \cap A_q$, moreover $\mathcal{C}_r(x) = \mathcal{C}_p(x)$ for $x \in A_p \setminus A_q$
and similarly $\mathcal{C}_r(x) = \mathcal{C}_q(x)$ for $x \in A_q \setminus A_p$ then $r = \<U_r,\mathcal{C}_r\>$ is
a common extension of $p$ and $q$ in $\mathbb{P}$.

To see that $r \in P$, the only non-trivial thing to check is that $\mathcal{C}_r(x) \subset ac_r(x)$ for any $x \in A_p \cap A_q$.
This, however, is clear because $U_r(x) = U_p(x) = U_q(x) \subset x+1$ and hence $C \in ac_r(x)$ iff $C \cap x \in ac_r(x)$,
moreover $\pi$ is the identity on $A_r \cap x = A_p \cap x =A_q \cap x$. Now, that $r$ extends both $p$ and $q$ is obvious.
\end{proof}

Using our assumption $\mathfrak{c} = 2^{<\mathfrak{c}}$, standard counting and delta-system arguments imply that
every subset $Q$ of $P$ with $|Q| = \mathfrak{c}^+$ contains two elements $p,q \in Q$ that satisfy the conditions
of Lemma \ref{lm:amalg1}, and hence are compatible. Consequently, $\mathbb{P}$ is indeed $\mathfrak{c}^+-CC$.

It is an immediate consequence of the above results that we have $\mathfrak{c}^V = \mathfrak{c}^{V^\mathbb{P}}$ and $(\mathfrak{c}^+)^V = (\mathfrak{c}^+)^{V^\mathbb{P}}$.

\begin{lemma}\label{lm:Aq1}
For every condition $q \in P$ and $\alpha \in \mathfrak{c}^+ \setminus \omega$  there is $p \le q$
such that $\alpha \in A_p$.
\end{lemma}

\begin{proof}[Proof of \ref{lm:Aq1}]
Indeed, if $\alpha \notin A_q$ then let $A_p = A_q \cup \{\alpha\}$ and extend $U_q$ to $A_p$
by putting $U_p(\alpha) = \omega \cup \{\alpha\}$. Also, we extend $\mathcal{C}_q$ by letting $\mathcal{C}_p(\alpha) = \emptyset$.
It is obvious then that $p = \<U_p,\mathcal{C}_p\>$ is as required.
\end{proof}

It immediately follows that if $G \subset P$ is $\mathbb{P}$-generic over $V$ then $U = \bigcup \{U_p : p \in G\}$ maps
$\mathfrak{c}^+$ into $[\mathfrak{c}^+]^\omega$ and the (obviously locally countable) topology $\tau$ generated by
the range of $U$ is crowded because
$$\mathcal{B} = \bigcup \{\mathcal{B}_p : p \in G\} \subs [\mathfrak{c}^+]^\omega$$
forms a base for $\tau$.

We still have to work to show that $\tau$ is SAU. The Hausdorff property of $\tau$ immediately follow from the
following result.

\begin{lemma}\label{lm:H1}
For every condition $q \in P$ and distinct $\alpha, \beta \in A_q$ there is $p \le q$ such that
for some $\gamma, \delta \in A_p$ we have $\alpha \in U_p(\gamma)$, $\,\beta \in U_p(\delta)$
and $U_p(\gamma) \cap U_p(\delta) = \emptyset$.
\end{lemma}

\begin{proof}[Proof of \ref{lm:H1}]
Start by fixing a countable $\tau_q$-open set $W$ containing both $\alpha$ and $\beta$, e.g.
$W = U_q(\alpha) \cup U_q(\beta)$ will work. For every $x \in W$ consider the following two subfamilies of $[W]^\omega$:
$$\mathcal{B}_x = \{B \in \mathcal{B}_q : x \in B \subset W\},\,\, \mathcal{A}_x = \{B \cap C : B \in \mathcal{B}_x \text{ and }  C \in \mathcal{C}_q(x)\}.$$
Then we have $|\mathcal{A}_x \cup \mathcal{B}_x| < \mathfrak{c}$ and hence for
$\mathcal{A} = \bigcup \{\mathcal{A}_x : x \in W\}$ and $\mathcal{B} = \bigcup \{\mathcal{B}_x : x \in W\}$
we have $|\mathcal{A} \cup \mathcal{B}| < \mathfrak{c}$ as well.

We may thus apply our assumption $\mathfrak{r} = \mathfrak{c}$ to the family $\mathcal{A} \cup \mathcal{B} \subset [W]^\omega$
to obtain a partition $W = E \cup F$ of $W$ such that $|E \cap A| = |F \cap A| = \omega$ for all $A \in \mathcal{A}$
and $|E \cap B| = |F \cap B| = \omega$ for all $B \in \mathcal{B}$.
Without loss of generality, we may also assume that $\alpha \in E$ and $\beta \in F$.

Let us now choose $\gamma, \delta \in \mathfrak{c}^+$ such that $A_q < \gamma < \delta$, and extend $U_q$
to $A_p = A_q \cup \{\gamma,\delta\}$ by putting $U_p(\gamma) = E \cup \{\gamma\}$ and $U_p(\delta) = F \cup \{\delta\}$.
It is clear from our choice of $E$ and $F$ that $\tau_p$ is crowded and $\mathcal{C}_q(x) \subs ac_p(x)$ for all $x \in A_q$.
Consequently, if we extend $\mathcal{C}_q(x)$ to $A_p$ by putting $\mathcal{C}_p(\gamma) = \mathcal{C}_p(\delta) = \emptyset$
then $p = \<U_p, \mathcal{C}_p\> \in P$ is as required.
\end{proof}

Now, it is trivial from Lemma \ref{lm:H1} that the generic topology $\tau$ in $V[G]$ is Hausdorff. Our next lemma will imply the only
missing thing required to show that $\tau$ is SAU, namely that any two infinite $\tau$-closed sets intersect.

\begin{lemma}\label{lm:SAU1}
Given any  condition $q \in P$ and two countably infinite sets $C,D \in [A_q]^\omega$, there is $p \le q$
such that both $C$ and $D$ belong to $\mathcal{C}_p(x)$ for some $x \in A_p$.
\end{lemma}

\begin{proof}[Proof of \ref{lm:SAU1}]
As in the previous proof, we start by fixing a countable $\tau_q$-open set $W$ such that $C \cup D \subs W$.
Then we choose $\gamma \in \mathfrak{c}^+$ such that $A_q < \gamma$ and extend $U_q$
to $A_p = A_q \cup \{\gamma\}$ by putting $U_p(\gamma) = W \cup \{\gamma\}$.
We also extend $\mathcal{C}_q(x)$
to $A_p$ by putting $\mathcal{C}_p(\gamma) = \{C, D\}$. It is trivial that then $p = \<U_p, \mathcal{C}_p\> \in P$
and $p \le q$, completing the proof.
\end{proof}

As a corollary we have that the generic topology $\tau$ is SAU.
Indeed, for any two sets $C, D \in [\mathfrak{c}^+]^\omega$, putting together Lemmas \ref{lm:c+cl1},  \ref{lm:Aq1}, and \ref{lm:SAU1}
we may conclude that
$$\{p \in P : C \cup D \subs A_p \text{ and } \exists\,x \in A_p \,\{C,D\} \subs ac_p(x)\}$$ is dense in $\mathbb{P}$.
But clearly if $\{C,D\} \subs ac_p(x)$ then $p$ forces $\{C,D\} \subs ac_\tau(x)$ and so
$cl_\tau(C) \cap cl_\tau(D) \ne \emptyset$ as well.
Thus the proof of Theorem \ref{tm:r=c} is completed.
\end{proof}

\medskip

We do not know if the above result is (consistently) sharp, i.e. if $\mathfrak{c}^+$ could be replaced by, say, $2^\mathfrak{c}$.
our next theorem shows that this is possible if we give up local countability. Also, the assumption $\mathfrak{r} = \mathfrak{c}$
in the ground model is strengthened to $\mathfrak{r}^* = \mathfrak{c}$, where $\mathfrak{r}^*$ is defined to be the smallest
cardinal $\varrho$ such that there is a family $\mathcal{A}$ of infinite sets that cannot be reaped (or bisected) by
a single set. So, $\mathfrak{r}^* = \mathfrak{c}$ just says that any family $\mathcal{A}$ of infinite sets
with $|\mathcal{A}| < \mathfrak{c}$ can be reaped. Clearly, we have $\omega < \mathfrak{r}^* \le \mathfrak{r} \le \mathfrak{c}$,
moreover $\mathfrak{r}^* < \mathfrak{r} = \mathfrak{c}$ is consistent.

\begin{theorem}\label{tm:r*=c}
If $\mathfrak{r^*} = \mathfrak{c} = 2^{<\mathfrak{c}}$ then there is a $\mathfrak{c}$-closed and
$\mathfrak{c}^+-CC$ notion of forcing $\mathbb{P}$ such that, in the generic extension $V^\mathbb{P}$, there is
a crowded SAU space of cardinality $(2^\mathfrak{c})^V = (2^\mathfrak{c})^{V^\mathbb{P}}$.
\end{theorem}

\begin{proof}
To simplify notation, we shall write $\kappa = 2^\mathfrak{c}$ in what follows.
Similarly as in the proof of Theorem \ref{tm:r=c}, our aim is then to force a function
$U : \kappa \to \mathcal{P}(\kappa)$ such that $\alpha \in U(\alpha)$ for all $\alpha \in  \kappa$
and the topology $\tau$ generated by the range $\{U(\alpha) : \alpha \in \kappa\}$ of $U$ on $\kappa$ is a crowded SAU space.
The notion of forcing $\mathbb{P} = \<P, \le \>$ will also be quite similar to the one used there.

The elements of $P$ will be pairs of the form
$p = \<U_p,\mathcal{C}_p\>$, where $U_p$ is a function with
domain $A_p \in [\mathfrak{c}^+]^{< \mathfrak{c}}$ with values taken in $\mathcal{P}(A_p)$ such that
$\alpha \in U(\alpha)$ for all $\alpha \in A_p$, moreover
the second co-ordinate $\mathcal{C}_p$ of the condition $p$ is also a function with domain $A_p$ but
such that $\mathcal{C}_p(x) \in [ac_p(x)]^{< \mathfrak{c}}$ for $x \in A_p$. Here we are freely using
the analogs of the pieces of notation from the above proof of \ref{tm:r=c}, so $ac_p(x)$ denotes the family of
all countable subsets of $A_p$ that accumulate to $x$ in the topology $\tau_p$ generated by the
range of $U_p$  on $A_p$. Also, we shall use the notation 
$B_{p,I} = \bigcap \{U_p(\alpha) : \alpha \in I\}$ for $I \in [A_p]^{<\omega}$
to obtain the base  $$\mathcal{B}_p = \{B_{p,I} : I \in [A_p]^{< \omega}\} \setminus \{\emptyset\}$$ for $\tau_p$.
 
Next we define the partial order $\le$ on $P$ that, as far as the first co-ordinate is concerned, is quite different from the corresponding
part of \ref{df:leq1}.

\begin{definition}\label{df:leq2}
For $p,q \in P$ we have $p \leq q$, i.e. $p$ is a stronger condition than $q$ iff
\begin{enumerate}[(a)]
\item $A_p \supset A_q$ and for every $\alpha \in A_q$ we have $A_q \cap U_p(\alpha) = U_q(\alpha)$;

\smallskip

\item $U_q(\alpha) \cap U_q(\beta) = \emptyset$ implies $U_p(\alpha) \cap U_p(\beta) = \emptyset$
for any $\alpha,\beta \in A_q;$

\smallskip

\item $\mathcal{C}_p(x) \supset \mathcal{C}_q(x)$ for all $x \in A_q$.
\end{enumerate}
\end{definition}
It is obvious that $\le$ is indeed a partial order on $P$.

We next present a series of lemmas that will help us prove the required properties of the forcing $\mathbb{P}$.

\begin{lemma}\label{lm:c+cl2}
The forcing $\mathbb{P}$ is $\mathfrak{c}$-closed.
\end{lemma}

\begin{proof}[Proof of \ref{lm:c+cl2}]
Assume that $\<p_\xi : \xi < \varrho\>$ is a decreasing $\varrho$-sequence in $\mathbb{P}$, where
$p_\xi = \< U_\xi, \mathcal{C}_\xi \>$ for $\xi < \varrho$ and $\varrho < \mathfrak{c}$  is an infinite regular cardinal.
We may then define a lower bound $p = \<U_p,\mathcal{C}_p\>$ for the $p_\xi$'s as follows.

\begin{enumerate}
\item  $A_p = \bigcup_{\xi < \varrho} A_\xi$, and $U_p(\alpha) = \bigcup \{U_\xi(\alpha) : \xi_\alpha \le \xi < \varrho\}$,
where $\xi_\alpha = \min\{\xi : \alpha \in A_\xi\}$.

\smallskip

\item $\mathcal{C}_p(x) = \bigcup \{\mathcal{C}_\xi(x) : x \in A_\xi \}$ for any $x \in A_p$.
\end{enumerate}

Since $\mathfrak{c}$  is regular, we have both $|A_p| < \mathfrak{c}$ and $|\mathcal{C}_p| < \mathfrak{c}$.
For any $\xi < \varrho$ and $C \in \mathcal{C}_\xi(x)$ we have $C \in ac_p(x)$
because then $C \in ac_\eta(x)$ for all $\eta \in \varrho \setminus \xi$ and every $B \in \mathcal{B}_p$
eventually belongs to $\mathcal{B}_\eta$ as well. Thus we have $p \in P$.
The easy verification of $p \le p_\xi$ for all $\xi < \varrho$ is left to the reader.
\end{proof}

As in the proof of \ref{lm:amalg1},
we next give an amalgamation result that will imply the $\mathfrak{c}^+-CC$ property of $\mathbb{P}$.

\begin{lemma}\label{lm:amalg2}
Assume that $p,q \in P$ are {\em isomorphic} conditions, i.e. $\,otp(A_p) = otp(A_q)$,
the unique order isomorphism $\pi : A_p \to A_q$ is the identity on $A = A_p \cap A_q$, moreover for all $\alpha \in A_p$
we have
\begin{enumerate}[(i)]
\item $U_q(\pi(\alpha)) = U_p(\alpha)$, and

\smallskip

\item $\mathcal{C}_q(\pi(\alpha)) = \{\pi[C] : C \in \mathcal{C}_p(\alpha)\}$.
\end{enumerate}
Then $p$ and $q$ are compatible in $\mathbb{P}$.
\end{lemma}

\begin{proof}[Proof of \ref{lm:amalg2}]
We define the desired common extension  $r = \<U_r,\mathcal{C}_r\>$ of $p$ and $q$ by the following stipulations:

\begin{enumerate}[(a)]
\item $U_r(x) = U_r(\pi(x)) = U_p(x) \cup U_q(\pi(x))$ for $x \in A_p$;

\smallskip

\item $\mathcal{C}_r(x) = \mathcal{C}_r(\pi(x)) = \mathcal{C}_p(x) \cup \mathcal{C}_q(\pi(x))$ for $x \in A_p$.
\end{enumerate}
As $\pi$  is the identity on $A$, both functions $U_r$ and $\mathcal{C}_r$ are well-defined on $A_r = A_p \cup A_q$.

To see that $r \in P$, we have to check that $\mathcal{C}_r(x) \subs ac_r(x)$ for all $x \in A_r$.
By symmetry, it suffices to do this for $x \in A_p$. In view of conditions (ii) and (b), what we have to show is that $x \in B \in \mathcal{B}_r$
implies $|B \cap C| = |B \cap \pi[C]| = \omega$ for any $C \in \mathcal{C}_p(x)$.
Now, any member of $\mathcal{B}_r$ is of the form $$B_{r, I \cup J} = \bigcap \{U_r(i) : i \in I\} \cap \bigcap \{U_r(j) : j \in J\}$$
with $I \in [A_p]^{< \omega}$ and $J \in [A_q]^{< \omega}$. It is easy to compute from 
condition (a) that
we have $$B_{r, I \cup J} = B_{p, I \cup \pi^{-1}[J]} \cup B_{q, \pi[I] \cup J},$$
moreover $B_{q, \pi[I] \cup J} = \pi[B_{p, I \cup \pi^{-1}[J]}]$.

\smallskip

Then $x \in B_{r, I \cup J}$ implies  $x \in B_{p, I \cup \pi^{-1}[J]}$, hence $|C \cap B_{p, I \cup \pi^{-1}[J]}| = \omega$,
consequently, then $|\pi[C] \cap B_{q, \pi[I] \cup J}| = \omega$ as well. But this clearly implies
$$|B_{r, I \cup J} \cap C| = |B_{r, I \cup J} \cap \pi[C]| = \omega,$$
just as we claimed.

To see that $r$ is a common extension of $p$ and $q$, by symmetry again, it suffices to show $r \le p$.
Clearly, only condition (b) of \ref{df:leq2} requires any checking for this. But as $\pi$ is an isomorphism,
$U_p(x) \cap U_p(y) = \emptyset$ implies $U_q(\pi(x)) \cap U_q(\pi(y)) = \emptyset$, moreover
$U_p(x) \cap A = U_q(\pi(x)) \cap A$ and $U_p(y) \cap A = U_q(\pi(y)) \cap A$, we indeed have
$U_r(x) \cap U_r(y) = \emptyset$.
\end{proof}

Since  $\mathfrak{c} = 2^{<\mathfrak{c}}$, standard counting and delta-system arguments imply that
every subset $Q$ of $P$ with $|Q| = \mathfrak{c}^+$ contains two isomorphic elements $p,q \in Q$ that satisfy the conditions
of Lemma \ref{lm:amalg2}, and hence are compatible. Consequently, $\mathbb{P}$ is indeed $\mathfrak{c}^+-CC$.

It is an immediate consequence of the above results that we have $\mathfrak{c}^V = \mathfrak{c}^{V^\mathbb{P}}$,
moreover by counting nice names we also may conclude $(2^\mathfrak{c})^V = \kappa = (2^\mathfrak{c})^{V^\mathbb{P}}$.

We also have the following immediate consequence of Lemma \ref{lm:amalg2} which will be used to show that the generic topology $\tau$
is (very) crowded.

\begin{lemma}\label{lm:crd}
For every $\alpha < \kappa$ the set of conditions
$$D_\alpha = \{r \in P : \forall\,B \in \mathcal{B}_r\,(B \setm \alpha \ne \emptyset)\}$$
is dense in $\mathbb{P}$.
\end{lemma}

\begin{proof}[Proof of \ref{lm:crd}]
Since $cf(\kappa) > \mathfrak{c}$, the domain $A_p$ of any condition $p \in P$ is bounded in $\kappa$.
So we can pick $q \in P$ isomorphic to $p$ such that $A_q > \{\alpha\} \cup A_p$. Now, let $r$
be the common extension of $p$ and $q$ constructed in Lemma \ref{lm:amalg2}. Then, as we have seen,
$$\mathcal{B}_r = \{B \cup \pi[B] : B \in \mathcal{B}_p\},$$
where $\pi$ is the unique order isomorphism from $A_p$ to $A_q$. Thus we clearly have $p \ge r \in D_\alpha$.
\end{proof}

That the generic map $U$ is defined on all of $\kappa$, follows from the following trivial lemma.

\begin{lemma}\label{lm:Aq2}
For every $\alpha < \kappa$ the set of conditions $E_\alpha = \{p \in P : \alpha \in A_p\}$ is dense in $\mathbb{P}$.
\end{lemma}

\begin{proof}[Proof of \ref{lm:Aq2}]
Indeed, if $q \in P$ and $\alpha \notin A_q$ then define $p = \<U_p,\mathcal{C}_p\>$ as follows:
Let $A_p = A_q \cup \{\alpha\}$, $\,U_p(x) = U_q(x)$ and $\mathcal{C}_p(x) = \mathcal{C}_q(x)$ for $x \in A_q$,
moreover $U_p(\alpha) = A_p$ and $\mathcal{C}_p(\alpha) = \emptyset$.
Then $p \le q$ and $p \in E_\alpha$.
\end{proof}

It immediately follows from Lemma \ref{lm:Aq2} that if $G$ is $\mathbb{P}$-generic over $V$ then putting
$U(\alpha) = \bigcup \{U_p(\alpha) : p \in G \cap E_\alpha\}$ for all $\alpha < \kappa$, then $U : \kappa \to \mathcal{P}(\kappa)$
is well-defined in $V[G]$, hence so is the topology $\tau$ generated by the range of $U$.
It is also clear that if $\mathcal{B}$ is the base of $\tau$ consisting of all non-empty finite intersections of the subbase
$\{U(\alpha) : \alpha < \kappa\}$ then every member of $\mathcal{B}$ is cofinal in $\kappa$, hence $\tau$ is crowded.
Indeed, fix $\alpha < \kappa$ and assume that $B_I = \bigcap \{U(i) : i \in I\} \ne \emptyset$ for some $I \in [\kappa]^{< \omega}$. By lemmas \ref{lm:Aq2}
and \ref{lm:crd} there is $p \in G \cap D_\alpha$ such that $I \subs A_p$ and $B_{p,I} \in \mathcal{B}_p$.
But then $B_{p,I} \subs B_I$ implies $B_I \setm \alpha \ne \emptyset$, so $ B_I$ is indeed cofinal in $\kappa$.

Our next result is used to show that the generic topology $\tau$ is Hausdorff.

\begin{lemma}\label{lm:H2}
For distinct $\alpha, \beta \in \kappa$ the set of conditions $D(\alpha,\beta)$ consisting of all $p \in P$ such that
for some $\gamma, \delta \in A_p$ we have $\alpha \in U_p(\gamma)$, $\,\beta \in U_p(\delta)$,
and $U_p(\gamma) \cap U_p(\delta) = \emptyset$ is dense in $\mathbb{P}$.
\end{lemma}

\begin{proof}[Proof of \ref{lm:H2}]
By Lemma \ref{lm:Aq2} it suffices to show that any $q \in P$ with $\alpha, \beta \in A_q$ has an extension in $D(\alpha,\beta)$.
For every $x \in A_q$ consider the following family of infinite sets:
$$ \mathcal{A}_x = \{B \cap C : x \in B \in \mathcal{B}_q \text{ and }  C \in \mathcal{C}_q(x)\} \subs [A_q]^\omega.$$
Then we have $|\mathcal{A}_x| < \mathfrak{c}$ and hence, as $\mathfrak{c}$ is regular,
$|\mathcal{A}| < \mathfrak{c}$ as well.

We may thus apply our assumption $\mathfrak{r}^* = \mathfrak{c}$ to $\mathcal{A} \subset [A_q]^\omega$
to obtain a partition $A_q = E \cup F$ of $A_q$ such that $|E \cap A| = |F \cap A| = \omega$ for all $A \in \mathcal{A}$.
Without loss of generality, we may also assume that $\alpha \in E$ and $\beta \in F$.

Let us now choose distinct $\gamma, \delta \in \kappa \setminus A_q$ and define $p \in P$ as follows.
Put $A_p = A_q \cup \{\gamma,\delta\}$, for $x \in A_q$ set $U_p(x) = U_q(x)$ and $\mathcal{C}_p(x) = \mathcal{C}_q(x)$,
moreover $U_p(\gamma) = E \cup \{\gamma\}$, $\,U_p(\delta) = F \cup \{\delta\}$ and, finally,  $\mathcal{C}_p(\gamma) = \mathcal{C}_p(\delta) = \emptyset$.
It is clear from our choice of $E$ and $F$ that $\mathcal{C}_p(x) \subs ac_p(x)$ for all $x \in A_q$. It trivially follows
then that $p = \<U_p, \mathcal{C}_p\> \in D(\alpha,\beta)$ and $p \le q$.
\end{proof}

It is obvious from condition (b) of Definition \ref{df:leq2} that every $p \in D(\alpha,\beta)$ forces $U(\alpha) \cap U(\beta) = \emptyset$,
hence $\tau$ is indeed Hausdorff. Our last result will finish our proof by implying that $\tau$ is SAU.

\begin{lemma}\label{lm:SAU2}
For any two sets $C, D \in [\kappa]^\omega$ the set of conditions
$$E(C,D) = \{p \in P : C \cup D \subs A_p \text{ and } \exists\,x \in A_p \,\{C,D\} \subs ac_p(x)\}$$ is dense in $\mathbb{P}$.
\end{lemma}

\begin{proof}[Proof of \ref{lm:SAU2}]
By Lemmas \ref{lm:c+cl2} and \ref{lm:Aq2} it suffices to show that any $q \in P$ with $C \cup D \subs A_q$ has an extension in $E(C,D)$.
To see that, pick $\gamma \in \kappa \setm A_q$ and define $p \in P$ as follows.
First, set $A_p = A_q \cup \{\gamma\}$, for $x \in A_q$ set $U_p(x) = U_q(x)$ and $\mathcal{C}_p(x) = \mathcal{C}_q(x)$.
We also set $U_p(\gamma) = \{\gamma\} \cup A_q$ and $\mathcal{C}_p(\gamma) = \{C, D\}$. Then the only
member of $\mathcal{B}_p$ containing $\gamma$ is $A_p$, hence we trivially have $\mathcal{C}_p(\gamma) \subs ac_p(\gamma)$.
We thus have $p = \<U_p, \mathcal{C}_p\> \in E(C, D)$ and $p \le q$, completing the proof.
\end{proof}

Since every condition $p \in E(C, D)$ clearly forces $cl_\tau(C) \cap cl_\tau(D) \ne \emptyset$,
it immediately follows that $\tau$ is indeed SAU.
\end{proof}

\bigskip

\section{Crowded SAU spaces from Cohen reals}

\bigskip

In the previous section, with considerable effort, we presented consistent examples of
crowded SAU spaces of size $> \mathfrak{c}$. In this section we show that to obtain such examples
of size $\le \mathfrak{c}$ is much easier, in fact, we can get them by
simply adding Cohen reals.

To fix our notation, we shall denote by $\mathbb{C}_I$ the standard Cohen forcing
$$\mathbb{C}_I = \<Fn(I,2), \supset\>,$$
using the notation of \cite{K1}.
Before turning to this promised result, we formulate and prove the following technical lemma.

\begin{lemma}\label{lm:cof}
Assume that $cof([\kappa]^\omega, \subset) = \kappa$, moreover $\tau$ is
a crowded and locally countable Hausdorff topology  on a set $X$ of cardinality $\kappa$ with $X \cap \kappa = \emptyset$.
Then there is a crowded and locally countable Hausdorff topology $\sigma$ on
$Z = X \cup \kappa$ in $V^{\mathbb{C}_\kappa}$ such that
\begin{enumerate}
\item $\tau \subs \sigma$;

\smallskip

\item $X \cap cl_\sigma(A) = cl_\tau(A)$ for every $A \in V \cap [X]^\omega$;

\smallskip

\item $cl_\sigma(C) \cap cl_\sigma(D) \ne \emptyset$ for any two $C, D \in V \cap [X]^\omega$.
\end{enumerate}
\end{lemma}

\begin{proof}\ref{lm:cof}
Since $X$ is locally countable, every countable subset of $X$ is included in a countable $\tau$-open set,
hence by $cof([\kappa]^\omega, \subset) = \kappa$ we can fix $\mathcal{U} = \{U_\alpha : \alpha \in \kappa\} \subs \tau \cap [X]^\omega$
that is cofinal in $[X]^\omega$. For each $\alpha \in \kappa$ we also fix an $\omega$-type one-one enumeration $\{x_{\alpha, n} : n < \omega\}$
of $U_\alpha$.

Let $G$ be $\mathbb{C}_\kappa$-generic over $V$ and $g = \cup G : \kappa \to 2$ be the corresponding Cohen generic map.
For every $\alpha < \kappa$ we then define $$W_\alpha = \{\alpha\} \cup \{x_{\alpha, n} : g(\omega \cdot \alpha + n) = 1\}.$$
We claim that the topology $\sigma$ generated in $V[G]$ on $Z$ by the subbase
$$\tau \cup  \{W_\alpha : \alpha < \kappa\} \cup \{T_\alpha = Z \setm W_\alpha : \alpha < \kappa\}$$
is as required.

That $\sigma$ is locally countable follows because $W_\alpha \subs U_\alpha$ for all $\alpha < \kappa$. To see that it is crowded,
consider any non-empty basic open set of the form $B = U \cap W_I \cap T_J$ where $|U| = \omega$, $\,W_I = \cap \{W_i : i \in I\}$
and $T_J = \cap \{T_j : j \in J\}$ with $I, J \in [\kappa]^{< \omega}$.  Let $H=U \cap \bigcap \{U_i:i\in I\}$.
Now $\empt\ne B\subs H$ implies that $H$ is infinite because $X$ is crowded.
This yields,  by genericity, that $|B| = \omega$ as well.

To check that $\sigma$ is Hausdorff, consider first $\alpha \ne \beta \in \kappa$. In this case $W_\alpha$ and $Z \setm W_\alpha$
are clearly their disjoint neighborhoods. If $x \in X$ and $\,\alpha \in \kappa$
then $|\{\beta : x \in U_\beta\}| = \kappa$ implies that for every condition $q \in Fn(\kappa, 2)$
there is $\beta \ne \alpha$ with $x = x_{\beta, n} \in U_\beta$ and $$[\omega \cdot \beta, \omega \cdot \beta + \omega) \cap \dom(q) = \emptyset.$$
Now, if $p$ extends $q$ so that $p(\omega \cdot \beta + n) = 1$ then $p$ forces $x \in W_\beta$ and $ \alpha \in Z \setm W_\beta$.
Finally, there is nothing to prove if $x \ne y \in X$ because
$\tau \subs \sigma$ is Hausdorff.

It remains to prove (1), (2) and (3). Now, (1) holds by definition. To check (2), consider any $x \in  cl_\tau(A)$ for some $A \in V \cap [X]^\omega$
and $x \in B = U \cap W_I \cap T_J$.
Let  $H = U \cap \bigcap \{U_i : i \in I \cup J\}$, then $x \in H \in \tau$.
We may assume that $x \notin A$, hence $|H \cap A| = \omega$.
So, given any condition $q \in Fn(\kappa, 2)$,  there is $y \in H \cap A$ such that
for every $i \in I \cup J$ with $y = x_{i,n_i}$ we have $\omega \cdot i + n_i \notin \dom(q)$.
Consequently, we may extend $q$ to $p \in Fn(\kappa, 2)$ so that $p(\omega \cdot i + n_i) = 1$ whenever $i \in I$
and $p(\omega \cdot i + n_i) = 0$ whenever $i \in J$. But then $p$ clearly forces $y \in B$, hence $B \cap A \ne \emptyset$ as well.

Finally, since for any $C, D \in V \cap [X]^\omega$ there is $\alpha < \kappa$ with $C \cup D \subs U_\alpha$,
to prove (3), it clearly suffices to show that $\alpha \in cl_\sigma(A)$ whenever $A \in V \cap [U_\alpha]^\omega$.
To see this, note first that the sets of the form $$B_{\alpha, I} = \{W_\alpha \setm \bigcup_{i \in I} W_i\}$$
constitute a $\sigma$-neighborhood base at $\alpha$, where $\alpha \notin I \in [\kappa]^{< \omega}$.
For any condition $q \in Fn(\kappa, 2)$  there is $x = x_{\alpha, n} \in A$ such that $\omega \cdot \alpha + n \notin \dom(q)$
and for every $i \in I$ if $x = x_{i,n_i} \in U_i$ then $\omega \cdot i + n_i \notin \dom(q)$.
So, we may extend $q$ to $p \in Fn(\kappa, 2)$ so that $p(\omega \cdot \alpha + n) = 1$ and $p(\omega \cdot i + n_i) = 0$
whenever $x \in U_i$ for some $i \in I$. But then $p$ forces $x \in B_{\alpha, I}$, hence $A \cap B_{\alpha, I} \ne \emptyset$,
completing the proof.
\end{proof}

From Lemma \ref{lm:cof} the following result is easily deduced.

\begin{theorem}\label{tm:Vi}
Assume that 
$W$ is a model of ZFC,  $\<V_{\alpha}:{\alpha}<{\omega}_1\>$ is an increasing sequence of ZFC submodels in $W$,
$\,(\omega_1)^{V_0} = (\omega_1)^W$, and
$\kappa \in V_0$ is a cardinal such that  $$({[{\kappa}]}^{{\omega}})^W=\bigcup_{{\alpha}<{\omega}_1}({[{\kappa}]}^{{\omega}})^{V_{\alpha}}$$ 
and        $cof([\kappa]^\omega, \subset) = \kappa$ holds in all the $V_\alpha$'s.
Assume also that $V_{\alpha+1}$ contains a $\mathbb{C}_{\kappa}$-generic filter $G_\alpha$ over $V_\alpha$ for every $\alpha < \omega_1$.
Then there is a crowded locally countable SAU space of cardinality $\kappa$ in $W$.
\end{theorem}

\begin{proof}
To start with, we fix a crowded locally countable space $\<X_0, \tau_0\> \in V_0$ of cardinality $\kappa$
such that $X_0 \cap (\omega_1 \times \kappa) = \emptyset$.
By transfinite recursion on $\alpha \le \omega_1$ we then define crowded locally countable topologies
$\tau_\alpha \in V_\alpha$ on $X_\alpha = \alpha \times \kappa$ as follows.

To obtain $\tau_{\alpha+1}$ from $\tau_\alpha$, we apply Lemma \ref{lm:cof} to get a crowded locally countable topology
$\sigma_\alpha \in V[G_\alpha] \subs V_{\alpha + 1}$ on $X_{\alpha+1} = X_\alpha \cup (\{\alpha\} \times \kappa)$
with properties (1) - (3) applied to $X_\alpha$ and $\{\alpha\} \times \kappa$ instead of $X$ and $\kappa$.
Then $\tau_{\alpha+1}$ is the topology generated by $\sigma_\alpha$ on $X_{\alpha+1}$ in $V_{\alpha + 1}$.
For $\alpha$ limit we simply let $\tau_\alpha$ be the topology on $X_\alpha$ generated by $\bigcup_{\beta < \alpha} \tau_\beta$
on $X_{\alpha}$ in $V_{\alpha +1}$.

Now, it is straight forward to check that $\tau_{\omega_1}$ is a crowded locally countable  topology on $X_{\omega_1}$
in the final model $W$. The topology is SAU because if $A,B\in {[X_{{\omega}_1}]}^{{\omega}}\cap W$
then there is ${\alpha}<{\omega}_1$ with $A,B\in {[X_{{\omega}_{\alpha}}]}^{{\omega}}\cap V_{\alpha}$, and so 
$A$ and $B$ have a common accumulation point in $X_{{\alpha}+1}$. Then 
$x$ is a common accumulation point in $X_{{\omega}_1}$ as well. 
\end{proof}

While it is an immediate corollary of Theorem \ref{tm:Vi} that the equality $cof([\kappa]^\omega, \subset) = \kappa$
implies the existence of a crowded locally countable
SAU space of cardinality $\kappa$ in the generic extension
obtained by adding $\kappa > \omega$ Cohen reals, it may be somewhat surprising that the two statements are actually equivalent.

\begin{theorem}\label{tm:cof}
For any uncountable cardinal $\kappa$ TFAE:
\begin{enumerate}[(i)]
\item $cof([\kappa]^\omega, \subset) = \kappa$.

\smallskip

\item There is a crowded locally countable SAU space in $V^{\mathbb{C}_\kappa}$.

\smallskip

\item There is a locally countable and countably compact $T_1$-space in some CCC generic extension $W$ of $V$.
\end{enumerate}
\end{theorem}

\begin{proof}
Since (i) $\,\Rightarrow\,$ (ii) is implied by Theorem \ref{tm:Vi} and (ii) $\,\Rightarrow\,$ (iii) is trivial,
it suffices to show that (iii) $\,\Rightarrow\,$ (i).

We first note that if $X$ is a locally countable and countably compact $T_1$-space of cardinality $\kappa > \omega$ then,
choosing a countable neighbourhood $U_x$ of every non-isolated point $x \in X'$, the family $$\mathcal{U} = \{U_x : x \in X'\} \subs [X]^\omega$$
has the property that for every $A \in [X]^\omega$ there is $U_x \in \mathcal{U}$ with $|A \cap U_x| = \omega$, i.e. $\mathcal{U}$ is $\omega$-hitting.

Now, it is well-known that if $\kappa > \omega$ then the existence of an $\omega$-hitting family $\mathcal{H}$ of size $\lambda$ in $[\kappa]^\omega$ implies that of
a cofinal family of size $\lambda$ in $[\kappa]^\omega$. Indeed, we may then take $\mathcal{H}$ with $|\mathcal{H}| = \lambda$
that is $\omega$-hitting in $[\kappa]^{< \omega}$ and then $\{\cup H : H \in \mathcal{H}\}$ is cofinal in $[\kappa]^\omega$.
This is because if $$A = \{\alpha_i : i < \omega\} \in [\kappa]^\omega \text{ \,and\, } S_A = \big\{\{\alpha_i : i < n\} : n < \omega\big\},$$
then $|H \cap S_A| = \omega$ implies $A \subs \cup H$.

Thus (iii) implies that $cof([\kappa]^\omega, \subset) = \kappa$ holds in $W$, but this, in turn, implies the same in $V$
because for any set $A \in W$ with $A \subs V$ there is $B \in V$ such that $A \subs B$ and $|A| = |B|$.
\end{proof}

It is easy to see that $\mathfrak{r} \ge \mathfrak{r}^* \ge \kappa$ holds in $V^{\mathbb{C}_\kappa}$, moreover Theorem 3.7. of \cite{JSSz}
as well as Theorem  \ref{tm:r=c} above used the assumption $\mathfrak{r} = \mathfrak{c}$ to obtain `"large" crowded locally countable SAU spaces.
This lead us to raise the question if one could get such spaces when $\mathfrak{r}$ is small. Our last result gives an affirmative answer to this question.

\begin{theorem}\label{tm:r=w1}
There are models of ZFC containing crowded locally countable SAU spaces of cardinality $\mathfrak{c}$ 
in which $\mathfrak{r} = \omega_1$ but $\mathfrak{c}$ is arbitrarily large.
\end{theorem}
 
\begin{proof}
To get such a model we first fix a cardinal $\kappa = \kappa^\omega$ in the ground model and then will do a
finite support iteration $\<\mathbb{P}_\alpha : \alpha \le \omega_1\>$ of length $\omega_1$ of CCC 
forcings 
where $\mathbb{P}_{\alpha+1} = \mathbb{C}_\kappa \ast \mathbb{Q}_\alpha$ for any $\alpha < \omega_1$.
Then, independently of the choice of the $\mathbb{Q}_\alpha$'s, we get from
 Theorem \ref{tm:Vi} that
a crowded locally countable SAU space of cardinality $\mathfrak{c} = \kappa$  exists in the final model $W = V^{\mathbb{P}_{\omega_1}}$.

The posets  $\mathbb{Q}_\alpha$ will be obtained together with ultrafilters $u_\alpha \in \omega^*$ in $V^{\mathbb{P}_\alpha}$ by recursion
so that $\beta < \alpha < \omega_1$ implies $u_\beta \subs u_\alpha$. Our $u_0$ is an arbitrary free ultrafilter on $\omega$ in the ground model and
for $\alpha$ limit we take $u_\alpha \supset \bigcup \{u_\beta : \beta < \alpha\}$. Once we have $u_\alpha$, 
we let $\mathbb{Q}_\alpha$ be the standard CCC, in fact even $\sigma$-centered, partial order that adds an
infinite pseudo-intersection $S_\alpha$ of $u_\alpha$.

Then $u = \bigcup \{u_\alpha : \alpha < \omega_1\}$ is a free ultrafilter on $\omega$ in the final model $W = V^{\mathbb{P}_{\omega_1}}$
that is generated by the family $$\{S_\alpha : \alpha < \omega_1\} \cup \{\omega \setm a : a \in [\omega]^{< \omega}\}.$$ 
Thus $W$ actually satisfies $\mathfrak{u} = \mathfrak{r} = \omega_1$,
see   \cite[V.4.27]{K1} for more details.
\end{proof}

\end{document}